\documentclass[10pt,a4paper]{amsart}
\usepackage{graphicx}
\usepackage{amssymb}
\usepackage{graphicx}
\usepackage[all]{xy}

\theoremstyle{plain}

\newtheorem{theorem}{Theorem}[section]
\newtheorem{corollary}[theorem]{Corollary}

\newtheorem{lemma}[theorem]{Lemma}
\newtheorem{proposition}[theorem]{Proposition}

\theoremstyle{definition}
\newtheorem{definition}[theorem]{Definition}
\newtheorem{example}[theorem]{Example}

\theoremstyle{remark}

\newtheorem{remark}[theorem]{Remark}
\newtheorem*{remark*}{Remark}

\numberwithin{equation}{section}
\numberwithin{figure}{section}

\newcommand{\cGa}{\Ga^\circ}
\newcommand{\ci}{i^\circ}

\newcommand{\op}[1]{\operatorname{#1}}
\newcommand\pd[2]{\frac{\pa #1}{\pa #2}}
\newcommand{\BL}{\biggl}
\newcommand{\BR}{\biggr}
\newcommand{\bl}{\bigl}
\newcommand{\br}{\bigr}

\newcommand{\Vect}{{\rm Vect}}
\newcommand{\abs}[1]{\lvert#1\rvert}
\newcommand{\norm}[1]{\left\| #1 \right\|}

\DeclareMathOperator\dist{dist}
\DeclareMathOperator\supp{supp}

\let\rom\textup

\def\CC{\mathbb{C}}
\def\RR{\mathbb{R}}
\def\lra{\longrightarrow}

\def\la{\lambda}
\def\e{\varepsilon}
\def\si{\sigma}
\def\ph{\varphi}
\def\om{\omega}
\def\Ga{\Gamma}
\def\Om{\Omega}
\def\La{\Lambda}
\def\cA{{\mathcal{A}}}
\def\cB{{\mathcal{B}}}
\def\cG{{\mathcal{G}}}
\def\cI{{\mathcal{I}}}
\def\cJ{{\mathcal{J}}}
\def\cK{{\mathcal{K}}}
\def\gS{{\mathfrak{S}}}
\def\ov{\overline}
\def\ovs{\overset}
\def\wt{\widetilde}
\def\wh{\widehat}
\def\pa{\partial}
\def\ciM{M^\circ}
\def\dvol{d{\op{vol}}}
\begin{document}

\title[Elliptic Theory on Manifolds with Corners: I]{Elliptic Theory
on Manifolds with Corners:\\ I. Dual Manifolds\\ and
Pseudodifferential Operators}

\author[Nazaikinskii]{Vladimir Nazaikinskii}

\address{%
Institute for Problems in Mechanics, Russian Academy of Sciences,
\newline\indent pr.~Vernadskogo 101-1, 119526 Moscow, Russia}

\email{nazaikinskii@yandex.ru}

\thanks{Supported in part by RFBR grants 05-01-00982
and 06-01-00098, by President of the Russian Federation grant
MK-1713.2005.1, and by the DFG project 436 RUS
113/849/0-1\protect\raisebox{1pt}{\circledR}\ ``K-theory and Noncommutative
Geometry of Stratified Manifolds"}
\author[Savin]{Anton Savin}

\address{Independent University of Moscow,
Bol'shoi Vlas'evskii per.~11,\newline\indent 119002 Moscow, Russia}

\email{antonsavin@mail.ru} \email{sternin@mail.ru}

\author[Sternin]{Boris Sternin}

\subjclass{Primary 58J40; Secondary 47G30, 32S05}

\keywords{Manifold with corners, elliptic operator, localization principle.}

\date{May 8, 2006}

\begin{abstract}
In this first part of the paper, we define a natural dual object for
manifolds with corners and show how pseudodifferential calculus on such
manifolds can be constructed in terms of the localization principle in
$C^*$-algebras. In the second part, these results will be applied to the
solution of Gelfand's problem on the homotopy classification of elliptic
operators for the case of manifolds with corners.
\end{abstract}
\maketitle
\tableofcontents
\section*{Introduction}

This paper deals with elliptic theory on manifolds with corners.

Such manifolds arise, e.g., if one supplements the class of closed manifolds by
manifolds with boundary and considers products of manifolds. A natural class of
operators on such manifolds was introduced by Melrose \cite{Mel1,Mel7}.
Operators on manifolds with corners have been actively studied,
e.g. see \cite{Bun2,Kra1,LaMo3,LeMo1,Loy1,MeNi2,MePi2,Mon3,MoNi1,Nis1}.

\medskip

The present paper consists of two parts. In the first part, we define a
natural dual object for manifolds with corners and show how
pseudodifferential calculus on such manifolds can be constructed in terms
of the localization principle in $C^*$-algebras. In the second part, these
results will be applied to the solution of Gelfand's problem on the
homotopy classification of elliptic operators for the case of manifolds
with corners.

In more detail, the outline of the first part is as follows. In
Sec.~\ref{par0}, we deal with the geometry of manifolds with corners.
Specifically,
\begin{itemize}
    \item In Sec.~\ref{par01} we
recall some facts and definitions concerning manifolds with corners. Most of
the material in this section is not new, except possibly in form.
    \item In Sec.~\ref{par02} we introduce a new geometric object, the dual
    manifold $M^\#$ of a manifold $M$ with corners, and study some structures on it.
    The importance of this space lies in the fact that, on the one hand,
    pseudodifferential operators on manifolds with corners can be naturally
    defined as operators local with respect to the action of the algebra of
    continuous functions on the dual manifold. On the other hand, as will
    be shown in the second part of paper, under an additional assumption
    the $K$-homology of the dual
    manifold $M^\#$ classifies the elliptic theory on $M$.
\end{itemize}

In Sec.~\ref{s2} we define zero-order pseudodifferential operators
($\psi$DO) in $L^2$ spaces on manifolds with corners. The definition is
based on the localization principle in $C^*$-algebras (e.g.,
see~\cite[Proposition~3.1]{PlSe6}), goes by induction over the depth of the
manifold, i.e., the maximum codimension of the strata (one starts from
smooth manifolds, which have depth zero), and naturally involves
parameter-dependent $\psi$DO (which serve as symbols for $\psi$DO at
subsequent inductive steps). Hence we need some preliminaries:
\begin{itemize}
    \item In Sec.~\ref{ss20} we introduce $L^2$ spaces on manifolds with
    corners.
    \item In Sec.~\ref{ss21} we discuss translation-invariant operators in
    vector bundles over manifolds with corners and their relationship with
    parameter-dependent operators.
    \item In Sec.~\ref{ss22} we present the
    adaptation~\cite{NaSaSt2} of the localization principle to operator families.
    The proofs are either contained in~\cite{NaSaSt2} or can be obtained from
    those in~\cite{NaSaSt2} by obvious modifications; hence we omit them altogether.
\end{itemize}
After that, in Sec.~\ref{ss23} we give the definition of $\psi$DO and prove
their properties.

\subsection*{Nomenclature}

We shall use the following notation.

$L^2(X,\mu,H)$ is the space of square integrable $H$-valued functions on a
metric space $X$ with respect to a measure $\mu$ (where $H$ is a Hilbert
space). We omit the argument $H$ if $H=\CC$ and also omit $\mu$ if it is
clear from the context.

$\cB H$ and $\cK H$ are the algebra of bounded operators and the ideal of
compact operators in a Hilbert space $H$.

$C(X,\cA)$ is the $C^*$-algebra of continuous bounded functions on $X$
ranging in the $C^*$-algebra $\cA$, and $C_0(X,\cA)$ is the subalgebra of
functions decaying at infinity. We omit the argument $\cA$ if $\cA=\CC$.

\section{Geometry}\label{par0}

\subsection{Manifolds with corners and their
faces\spacefactor1001}\label{par01}

\begin{definition}
A \textit{manifold of dimension $n$ with corners} is a Hausdorff topological
space $M$ in which each point $x$ has a coordinate neighborhood of the form
$\ov\RR{}_+^d\times\RR^{n-d}$, $d=d(x)\in\{0,\dotsc,n\}$ where $x$ is
represented by the origin. Moreover, the transition maps are smooth functions.
Unless specified otherwise, we assume that $M$ is connected and compact. The
maximum number $d$ is called the \emph{depth} of the manifold and will be
denoted by $k=k(M)$.
\end{definition}

\begin{figure}
\begin{center}
\includegraphics[height=6cm]{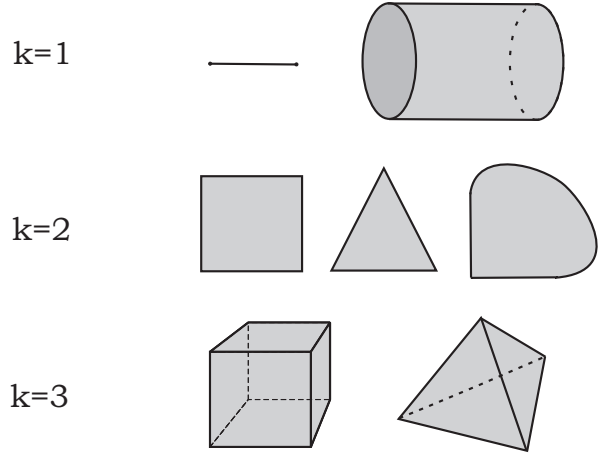}
\caption{Manifolds with corners of depth $k$.}\label{ris1}
\end{center}
\end{figure}

Some examples of manifolds with corners are shown in Fig.~\ref{ris1}.

\subsubsection*{Local defining functions}

By definition, each point $x\in M_j$ has a neighborhood $U\subset M$ with
local coordinates $\rho_1,\ldots,\rho_n$ such that the manifold is
determined in these coordinates by the system of inequalities
\begin{equation}\label{ineq1}
\rho_1\ge 0,\ldots,\rho_j\ge 0.
\end{equation}
The coordinates $(\rho_1,\ldots,\rho_j)$ are called \emph{defining functions}
of the face $\overline{\mathbb{R}}^j_+\subset \overline{\mathbb{R}}^n_+$.

\subsubsection*{Faces}

The set
\begin{equation*}
 \{x\in M\colon d(x)=l\}
\end{equation*}
is a smooth manifold of codimension $l$ in $M$. Its connected components are
called \textit{open faces of codimension} $l$. Let $\cGa_j(M)$, $j=0,\dotsc,N$
be all possible open faces of $M$, and let $d_j$ be their codimensions. We
assume that $d_0=0$ (so that $\cGa_0(M)=\ciM$ is the interior of $M$) and
$d_j>0$ for $j>0$. Thus $M$ is represented as the disjoint union
\begin{multline*}
 M=\bigsqcup_{j=0}^N\cGa_j(M)\equiv\ciM\sqcup\pa M,\\
\text{ where
    $\displaystyle\pa M=\bigsqcup_{j=1}^N\cGa_j(M)$ is the \textit{boundary}
    of
    $M$}.
\end{multline*}
Faces of codimension one are called \emph{hyperfaces}.

\begin{proposition}\label{p-1}
There exist canonically defined manifolds $\Ga_j(M)$ with corners such that
$\cGa_j(M)$ is the interior of $\Ga_j(M)$ and the diagram
\begin{equation*}
    \xymatrix{
  \cGa_j(M) \ar[dr]_{\ci_j} \ar[r]
                & \Ga_j(M) \ar@{-->}[d]^{i_j}  \\
                & M,             }
\end{equation*}
where the horizontal arrow and $\ci_j$ are natural embeddings and $i_j$ is
an immersion of manifolds with corner, commutes. The manifold $\Ga_j(M)$ is
called a \textit{closed face} of $M$.
\end{proposition}

\begin{proof}
[Proof \/\rm is given in the Appendix.]
\end{proof}

Since $\Ga_j(M)$ is a compact manifold with corners, we have
\begin{equation*}
 \pa\Ga_j(M)=\bigsqcup_{l=1}^L\cGa_l(\Ga_j(M)).
\end{equation*}
The image under $i_j$ of each open face $\cGa_l(\Ga_j(M))$, $l>0$, of th
manifold $\Ga_j(M)$ coincides with some open face $\cGa_r(M)$, $r=r(l)$, of
$M$ with $d_r>d_j$. in this case, we say that the face $\Ga_j$ (or
$\cGa_j$) and $\Ga_r$ (or $\cGa_r$) are \emph{adjacent to each other} and
write $\Ga_j\succ\Ga_r$. The restriction
\begin{equation*}
    i_{jl}:=i_j|_{\cGa_l(\Ga_j(M))}\colon \cGa_l(\Ga_j(M))\lra \cGa_{r(l)}(M)
\end{equation*}
is a finite covering whose structure group is a quotient of the homotopy
group $\pi_1(\cGa_{r(l)}(M))$ and a subgroup of the permutation group
$\gS_m$, where $m$ is the number of sheets of the cover.

\subsubsection*{The compressed cotangent bundle}

The compressed cotangent bundle $T^*M$ of a manifold $M$ with corners is
defined in the usual way (see~\cite{Mel1}). We take the subspace
$\Vect_b(M)$ of the space $\Vect(M)$ of vector fields on $M$ formed by
vector fields tangent to all open faces. The subspace $\Vect_b(M)$ is a
locally free $C^\infty(M)$-module.

Indeed, in local coordinates
$$
(\rho_1,\dotsc,\rho_d,y_{d+1},\dotsc,y_n)\in\ov\RR{}_+^d\times\RR^{n-d}
$$
a local basis in $\Vect_b(M)$ is formed by the vector fields
\begin{equation*}
    \rho_1\pd{}{\rho_1},\dotsc,\rho_d\pd{}{\rho_d},
    \pd{}{y_{d+1}},\dotsc,\pd{}{y_{n}}.
\end{equation*}

Consequently, $\Vect_b(M)$ is the section space of some vector bundle on
$M$, which will be denoted by $TM$ (the \textit{extended cotangent bundle}
of $M$), and the \emph{compressed cotangent bundle} $T^*M$ is now defined
as the bundle ($\RR$-) dual to $TM$. In the local coordinates
$(\rho_1,\dotsc,\rho_d,y_{d+1},\dotsc,y_n)$, a basis in the module of
sections of $T^*M$ is given by the forms
\begin{equation*}
    \rho_1^{-1}d\rho_1,\dotsc,\rho_d^{-1}d\rho_d,dy_{d+1},\dotsc,dy_n.
\end{equation*}

\subsubsection*{Conormal bundles of faces}

Let $ F=\cGa_j(M)$ be an open face of codimension $d=d_j$ in $M$. We define the
\textit{conormal bundle} of $ F$ as the subset $N^* F\subset T^*M|_{ F}$ formed
by functionals $\xi$ vanishing on any vector $v\in TM|_{ F}$ that can be
continued to a vector field second-order tangent to all faces in $\pa M$. One
readily sees that $N^* F$ is indeed a vector bundle; a basis in its fiber
consists of the $1$-forms
\begin{equation*}
    \rho_1^{-1}d\rho_1,\dotsc,\rho_d^{-1}d\rho_d.
\end{equation*}
This bundle can be canonically extended to a bundle over the closed face
$\ov{F}$; the latter bundle is called the \emph{conormal bundle} of
$\ov{F}$ and is denoted by $N^*\ov{F}$.

\begin{proposition}
One has the canonical direct sum decomposition
\begin{equation*}
    T^*M|_{\ov{F}}=T^*\ov{F}\oplus N^*\ov{F}
\end{equation*}
\rom{(}where the bundle on the left-hand side is obtained as the pullback
under the immersion of  $\ov{F}$ in $M$\rom{)}.
\end{proposition}
\begin{proof}
The assertion is local, so that we can assume that $\ov{F}$ is embedded in
$M$. Then the embedding  $T^*\ov{F}\subset T^*M$ is obtained as the map
dual to the restriction
$$
\Vect_b(M)\lra\Vect_b(\ov{F})
$$
of vector fields in $\Vect_b(M) $to $\ov{F}$. Now the desired properties
can be verified in coordinates.
\end{proof}

\subsubsection*{Normal bundles of faces}

Let $ F=\cGa_j(M)$ be again an open face of codimension $d=d_j$ in $M$, and
let $(\rho,y)$ and $(\wt\rho,\wt y)$ be two coordinate systems on $M$ in a
neighborhood of some point in $ F$. Since $\rho=\wt\rho=0$ on $ F$, we see
that the change of variables $(\rho,y)\longmapsto(\wt\rho,\wt y)$ has the
form
\begin{equation}\label{chu}
    \wt y=f(y)+O(\rho),\quad \wt\rho=A(y)\rho+O(\rho^2),
\end{equation}
where $A(y)$ is a smooth $d\times d$ matrix function. The
mapping~\eqref{chu} should take the positive quadrant with respect to the
variable $\rho$ to itself, and hence, letting $\rho$ tend to zero, we
verify that
\begin{equation*}
    A(y)=\Pi(y)\La(y),
\end{equation*}
where $\Pi(y)$ is a permutation matrix (and hence is locally constant in
$y$) and
$$
\La(y)=\op{diag}\{\la_1(y),\dotsc,\la_d(y)\}
$$
is a diagonal matrix with positive entries. The cocycle condition for the
matrices $A(y)$ implies that the matrices $\Pi(y)$ themselves satisfy the same
cocycle condition, so that we can define the $d$-dimensional real vector bundle
$N F$ over $ F$ for which the matrices $\Pi(y)$ are the transition mappings.
The change of variables
$$
t_j=-\ln\rho_j,\quad j=1,\dots,d,
$$
clarifies the meaning of this bundle. The second component in~\eqref{chu}
becomes
\begin{multline*}
    \wt t=\Pi(y) t+\ln\La(y)+O(e^{-2t})=\Pi(y) t+O(1), \\
    t_j\to+\infty,\quad j=1,\dots,d.
\end{multline*}
Thus $N F$ is the ``bundle of logarithms of determining functions'' of the
submanifold $ F$. We call it the \textit{logarithmic normal bundle} of $
F$. The matrices $\Pi$ simultaneously specify a bundle of positive
quadrants $\overline{\RR}_+^d$ over $ F$, which we denote by $N_+ F$ and
call the \textit{normal bundle} of $ F$. We have the exponential mapping
\begin{align*}
    \exp:N F &\lra N_+ F ,\\
    (y,t)&\longmapsto (y,\exp(-t))=(y,e^{-t_1},\dotsc,e^{-t_d}),
\end{align*}
which diffeomorphically maps the first bundle onto the interior of the
second.

One readily sees that both bundles naturally extend to bundles  $N\ov{F}$
and $N_+\ov{F}$ over the closed face $\ov{F}$.

By construction, the structure group of these bundles is a subgroup
$\gS_{\ov{F}}$ of the permutation group $\gS_d$. (Thus the numbering of the
coordinates $\rho$ in all charts is chosen in such a way that the
transition matrices range in the subgroup $\gS_{\ov{F}}$).
\begin{remark*}
We shall assume that the bundles $N\ov{F}$ and $N^*\ov{F}$ are reduced to
the minimal possible permutation structure group  $\gS_{\ov{F}}$. This will
be used in the sequel (in particular, see Lemma~\ref{eq2}).
\end{remark*}

\begin{proposition}\label{hahaha}
The logarithmic normal bundle $N\ov{F}$ and the conormal bundle $N^*\ov{F}$
are canonically dual.
\end{proposition}
\begin{proof}
[Proof \/\rm will be given in the Appendix]
\end{proof}

\begin{remark*}
The bundles $N\ov{F}$ and $N^*\ov{F}$ viewed as bundles with the structure
group $\gS_{\ov{F}}$ are canonically isomorphic, since permutation matrices
are unitary.
\end{remark*}

\subsubsection*{Compatible exponential mappings}

For each closed face $\Ga_j(M)$ of a manifold $M$ with corners, we have
defined the normal bundle $N_+\Ga_j(M)$. Just as with submanifolds of
smooth manifolds, one can define \textit{exponential mappings} of these
bundles into the manifold $M$ itself, which are local diffeomorphisms in a
neighborhood of the zero section. Moreover, for adjacent faces these
diffeomorphisms will be compatible in some sense. More precisely, the
following theorem holds.
\begin{theorem}\label{papkahoroshaja}
Let $\e>0$ b sufficiently small. Then there exist smooth mappings
\begin{equation*}
    f_j\colon N_+\Ga_j(M)\lra M,\quad j=1,\dotsc,N,
\end{equation*}
defined for $\abs{\rho}\le\e$, where $\rho$ is the variable in the fiber of the
bundle $N_+\Ga_j(M)$, such that the following conditions hold\rom:
\begin{enumerate}
    \item[\rom1.] On the zero section, $f_j=i_j$.
    \item[\rom2.] $f_j$ is a local diffeomorphism.
    \item[\rom3.] The restriction $f_j|_U$ of the mapping $f_j$ to
    some neighborhood of the open face $\cGa_j(M)$ in $N_+\Ga_j(M)$
    is a diffeomorphism.
    \item[\rom4.] If $\Ga_j(M)\succ\Ga_l(M)$, then the mappings $f_j$ and $f_l$
    are locally compatible in the following sense. In a neighborhood of
    any point
    $x\in\Ga_l(M)$, the diagram
\begin{equation}\label{dia1}
     \xymatrix{
  N_+\Ga_j(M) \ar[d]_{\pi_1} \ar[r]^{f_l^{-1}\circ f_j}
                & N_+\Ga_l(M) \ar[d]^{\pi_2}  \\
  \Ga_j(M) \ar[r]_{f_l^{-1}\circ f_j}
                & \ph(\Ga_j(M))             }
\end{equation}
commutes, where $\pi_1$ is the natural projection and $\pi_2$ is the
projection in the fibers of $N_+\Ga_l(M)$ onto the coordinate subbundle
into which $\Ga_j(M)$ is mapped under the local diffeomorphism
$\varphi=f_l^{-1}\circ f_j$, along the complementary coordinate subbundle.
\end{enumerate}
\end{theorem}

\begin{proof}
[Proof \/\rm will be given in the Appendix]
\end{proof}

\begin{remark*}
(a) Let $\Ga_j(M)\succ\Ga_r(M)$. Since
\begin{equation*}
     N_+\cGa_l(\Ga_j(M))\subset i_{jl}^{*}N_+\Ga_r(M),\quad r=r(l),
\end{equation*}
we see that by specifying a compatible tuple of exponential mappings $f_j$
for the strata of $M$ we automatically specify such tuples for the strata
of any closed stratum of $M$.

(b) the composition
$$
\wt f_j=f_j\circ\{t\mapsto e^{-t}\}\colon N\Ga_j(M)\lra M
$$
will also be referred to as the exponential map.
\end{remark*}

\begin{corollary}\label{coords1}
The manifold $M$ can be covered by finitely many coordinate neighborhoods
$U$ with coordinates $\rho_U=(\rho_1,\ldots,\rho_n)$ such that $M$ is given
in these coordinates by the system of inequalities~\eqref{ineq1} and the
following compatibility condition holds. Suppose that two charts $U$ and
$U'$ have a nonempty intersection.
\begin{enumerate}
    \item If the number of defining functions in $U$ and $U'$ is the same, then they
    coincide
    in $U\cap U'$ up to a permutation;
    \item otherwise, the smaller set of defining functions is a subset of
    the larger set in $U\cap U'$.
\end{enumerate}
\end{corollary}

\begin{remark}
This assertion plays in the theory of manifolds with corners the same role
as the collar neighborhood theorem in the theory of manifolds with boundary
(and contains the latter for the case in which the depth $k$ is equal to
one.
\end{remark}

To be definite, we assume in the following that the defining functions
specify coordinates in the domain where they are less than $3/2$.

\subsection{The dual manifold $M^\#$ and the algebra
$C(M^\#)$\spacefactor1001}\label{par02}

\subsubsection*{Definitions}

On the space $\RR_t^k\times\RR_x^m$, we define the algebra $C(k,m)$ of
bounded continuous functions $f(t,x)$ such that
\begin{equation*}
    f(\om\abs{t},x)\lra F(\om)\quad\text{as $\abs{t}\to\infty$}
\end{equation*}
uniformly with respect to $x$ and $\om=t/\abs{t}$, where $F(\om)$ is some
(continuous) function.

In the algebra of continuous functions on the interior $\ciM$ of the
manifold $M$, we single out a subalgebra $C(M^\#)$ as follows. We say that
$f\in C(M^\#)$ if for each coordinate neighborhood $U\simeq
\RR_+^k\times\RR^{n-k}$ on $M$ the function
\begin{equation*}
    F(t,x)=f|_U(e^{-t_1},\dotsc,e^{-t_k},x_{k+1},\dotsc,x_{n})
\end{equation*}
can be extended to a function in $C(k,n-k)$.

One can readily see that each function $f\in C(M^\#)$ is constant on each
hyperface of $M$.

One can readily describe the space $M^\#$ of maximal ideals of the algebra
$C(M^\#)$. As a set, it is the disjoint union of the interior $\ciM$ of the
manifold $M$ and the following sets $F^\#$ corresponding to faces $F$ of
positive codimension.
\begin{itemize}
    \item To each hyperface $F$, there corresponds a singleton $F^\#$.
    \item To each face $F$ of codimension $k=\op{codim} F>1$, there corresponds a
    set $F^\#$ that is the quotient of the open $k-1$-simplex
\begin{equation*}
    \ovs{\circ}\triangle_{k-1}=\left\{x\in\RR^k\colon x_i>0,\,i=1,\dotsc,k,\quad
    \sum x_i=1\right\}
\end{equation*}
    by the action of the structure group $\gS_{F}$ of the bundle
    $NF$.
\end{itemize}

The topology on $M^\#$ can be defined as follows. A sequence $z_n\in\ciM$
converges to a point $z\in F^\#$ if for each $\e>0$  all points $z_n$ lie
in the image of the neighborhood $\{\abs{\rho}<\e\}\subset N_+F$ starting
from some moment and
$$
\dist(\ln r_1(z_n):\dotsm:\ln r_k(z_n),z)\to0.
$$
(In the last formula, the line, treated as a point of the projective space,
is identified with the point of the simplex through which it passes, and
$z$ is understood as a $\gS_F$-orbit, $z\subset
\ovs{\circ}\triangle_{k-1}$.) Finally, the adjacency conditions for the
sets $F^\#$ naturally follow from the adjacency of the corresponding faces
and are induced by embeddings of simplices of various dimensions.

\begin{example}
(1) If $M$ is a manifold with boundary, then $C(M^\#)$ is the algebra of
functions constant on connected components of the boundary. Hence $M^\#$ is
obtained by retracting each boundary component into a point.

(2) Duality relates the cube to the octahedron and the icosahedron to the
dodecahedron. The tetrahedron proves to be self-dual.
\end{example}

\subsubsection*{Fibered structure on the dual space}\label{fiber}

Here we assume that $M$ is a manifold with corners such that the normal
bundle of each face is trivial and show that a neighborhood of each simplex
$F^\#$ of the dual manifold ${M}^\#$ is fibered over $F^\#$ with fiber
being a cone. This result will be used only in the proof of the
classification theorem in the second part of this paper.

Let $F\subset M$ be an open face of codimension $j$. We shall construct a
neighborhood ${U}^\#$ of the simplex ${F}^\#$ in the dual manifold
${M}^\#$.

First, we construct a neighborhood $U\subset N_+F$. It is convenient to use
the logarithmic coordinates
\begin{equation}\label{loga1}
\begin{array}{rcl}
  \ln: N_+F\setminus F & \stackrel\simeq\longrightarrow & NF, \\
   (x,\rho_F)& \mapsto  & (x,y=-\ln \rho_F). \\
\end{array}
\end{equation}
Here $\rho_F=(\rho_1,\ldots,\rho_j)$ is the set of defining functions of
$F$. By virtue of the triviality assumption, it is globally defined.

The image of the set in which $\rho_l<1$ for all $1\le l\le j$ will be
denoted by $N'_+F\subset NF$. In the coordinates $y$, it is given by the
condition $y>0$.

We use similar coordinates in neighborhoods of all faces of the face
$\overline{F}$. Then in the space $N'_+F$ we obtain the following
coordinates: the coordinates $y\in {\mathbb{R}}^j_+$ in the fibers; the
coordinates in the neighborhood $\mathbb{R}_+^{l}\times
\mathbb{R}^{n-j-l}\subset F$ of codimension $l$ in $F$, which will be
denoted by
$$
(x,\omega), \text{ where }(x_1,\ldots,x_l)=-\ln(\rho_{1},\ldots,\rho_l).
$$
(The coordinates $x$ are uniquely determined up to permutation; the number of
these coordinates is determined by the codimension in $F$ of the face near
which the point sits.)

To construct the neighborhood $U$, on ${F}$ we define the function
$|x|:=\sum_s x_s$. This is invariant under permutations of defining
functions and hence well defined.

\begin{figure}
\begin{center}
\includegraphics[height=9cm]{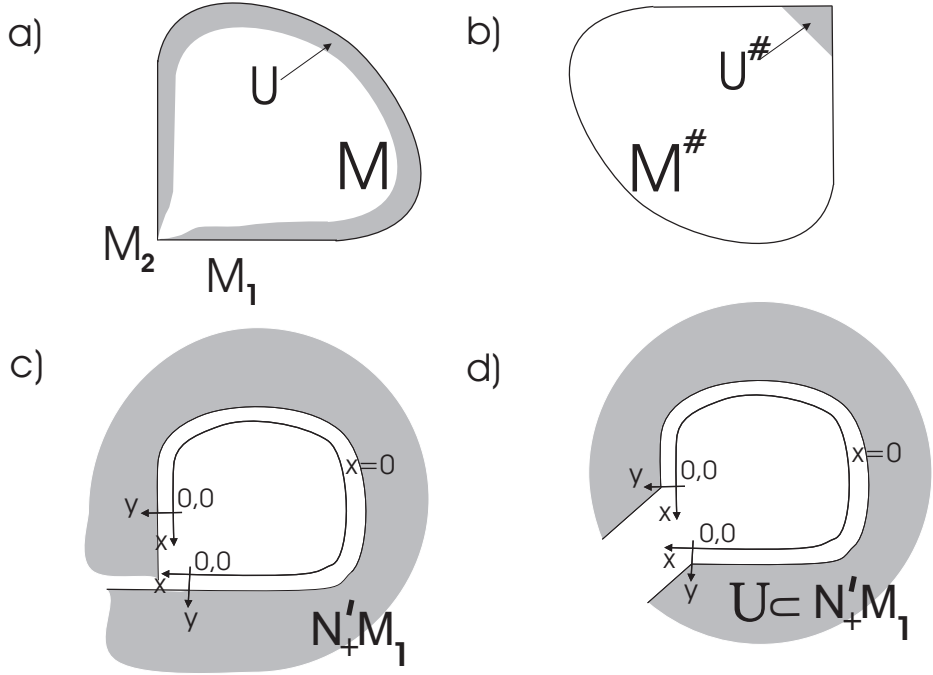}
\caption{a) the manifold $M$; b) the dual space
${M}^\#$; c) the positive quadrant in the normal bundle $N'_+M_1$; d) the
neighborhood $U\subset N'_+M_1$.}\label{r1}
\end{center}
\end{figure}

Now we globally define a set $U\subset N_+F$ by the condition
$$
U=\Bigl\{(y,x,\omega)\in N'_+F \;|\; \min y> |x|+1\Bigr\}
$$
in local coordinates, where $\min y$ is the minimum of the coordinates
$y_1,\ldots,y_j$. By way of example, Fig.~\ref{r1} shows the case in which
the manifold with corners is a $1$-gon; the set $U$ corresponding to the
one-dimensional edge is shown in the lower part of the figure as a dashed
infinite domain.

Consider the space
$$
{M}^\#_{\ge j}={M}^\#\setminus \bigcup^{j-1}_{j'=1}{M}^\#_{j'},
$$
obtained from ${M}^\#$ by deleting all simplices of dimension $\le j-2$.
\begin{lemma}
\begin{enumerate}
    \item the restriction of the projection $p:N_+F\to M$ to $U$ is one-to-one
    \rom{(}i.e.,
$U$ can also be treated as an open set in $M$; see top left in
Fig.~\rom{\ref{r1})};
    \item The dual space ${U}^\#\subset {M}^\#$ is an open neighborhood
    of the open simplex ${F}^\#$ in ${M}^\#_{\ge j}$
    \rom{(}see top right in Fig.~\rom{\ref{r1})}.
\end{enumerate}
\end{lemma}
\begin{proof}
Let us prove that $p$ is one-to-one. This can be violated only where
distinct parts of $F$ meet each other. We should prove that the projections of
components of $U$ corresponding to two adjacent faces are disjoint. Indeed,
let $U$ be defined in the first part by the condition
$$
\min y> |x|+1.
$$
Then in the second part some of the coordinates $x_I$ are interchanged with
some of the coordinates $y_I$ for some nonempty index set $I$. Then the set
$U$ in the second part is described by the inequality
$$
\min (x_I,y_{\overline{I}})>|y_I|+|x_{\overline{I}}|+1
$$
(in the original coordinates). Writing out these two systems componentwise,
we see that they are inconsistent, so that the projections of the parts of
$U$ into $M$ are disjoint.

The second assertion holds by construction.
\end{proof}

Now we can prove that the neighborhood ${U}^\#$ of the stratum ${F}^\#$ is
homeomorphic to the product of ${F}^\#$ by the cone
$$
K_{\Omega}=[0,1)\times \Omega/\{0\}\times \Omega.
$$
Here the base $\Omega$ of the cone is the dual space ${\overline{F}}^\#$ of
the closed face $\overline{F}$. The dual manifold is well defined, since
the closed face is a manifold with corners. As a result, we find that
${M}^\#$ is a stratified manifold with singularities.

\begin{proposition}\label{fiber1}
The projection
$$
\widetilde{p}:{U}^\#\to {F}^\#,\quad  \widetilde{p}(y):= y/|y|,
$$
is well defined. Its fiber is the cone $K_{{{\overline{F}}^\#}}$, and there
is a homeomorphism
$$
{U}^\#\simeq {F}^\#\times K_{{\overline{F}}^\#}.
$$
\end{proposition}
\begin{proof}
A straightforward computation shows that the projection $\widetilde{p}$ is
well defined. We define the map
$$
\begin{array}{ccc}
U &\longrightarrow & F^\#\times (0,1)\times F,\\
(y,x,\omega) &\mapsto & \left[\displaystyle\frac y{|y|}, \frac{|x|+1}{\min
y},(x,\omega)\right],
\end{array}
$$
The inverse map has the form
$$
\begin{array}{ccc}
F^\#\times (0,1)\times F &\longrightarrow & U,\\
(\theta,r,x,\omega) &\mapsto &
\left[\displaystyle\frac\theta{\min\theta}\frac{|x|+1}r,x,\omega\right],
\end{array}
$$
A routine verification of the fact that these mappings extend to
homeomorphisms ${U}^\#\simeq F^\#\times K_{{\overline{F}}^\#}$ is left to
the reader.
\end{proof}

\section{Pseudodifferential Operators}\label{s2}

\subsection{The space $L^2(M)$}\label{ss20}

Our $\psi$DO will act in the space $L^2(M)$, which is defined as follows.

Let $M$ be a compact manifold with corners, and let $\dvol$ be a smooth
measure on $M$ (obtained, say, via  an embedding of $M$ in a compact
Riemannian manifold). Now for each point $x\in M$ we define a measure
$\mu_x$ in $\ciM\cap U_x$, where $U_x\simeq
V\subset\ov\RR_+^k\times\RR^{n-k}$, $k=d(x)$, is a coordinate neighborhood
of $x$, by setting
\begin{equation}\label{mera-loc}
    \mu_x=(\rho_1\rho_2\dotsm \rho_k)^{-1}\dvol,
\end{equation}
where $\rho_1,\dotsc,\rho_k$ are the coordinates in the $\RR_+$-factors.
Next, we take a finite cover $M=\bigcup_{j=1}^{N'}U_{x_j}$ and a
subordinate partition of unity $\{e_j\}$ and set
\begin{equation}\label{mera}
    \mu=\sum_{j=1}^{N'}e_j\mu_{x_j}.
\end{equation}
This measure is up to equivalence independent of the ambiguity in the
construction, and hence the space $L^2(M)\ovs{\op{def}}\equiv
L^2(\ciM,\mu)$ is well defined up to norm equivalence. For the following,
we choose and fix such a measure and hence a Hilbert space structure in
$L^2(M)$. Note that the interiors of $M$ and $M^\#$ are the same, and so
$L^2(M)$ can also be viewed as $L^2(M^\#)$ (with respect to the same
measure). Hence it bears the natural structure of a $C(M^\#)$-module.

\subsection{Translation-invariant operators}\label{ss21}

In this section we shell work in the category whose objects are
arithmetic spaces $\RR^s$ and whose
morphisms are linear mappings taking each standard basis vector to zero or
to some standard basis vector. In particular, the automorphism group of
$\RR^s$ in this category is exactly the permutation group $\gS_s$ on the
$s$ standard basis vectors.

Let $M$ be a connected compact manifold with corners, and let $E\lra M$ be
a bundle with fiber $\RR^s$ on $M$. We reduce $E$ to a minimal structure
group $\cG\subset\gS_s$ (which is uniquely determined up to conjugacy) and
consider the principal $\cG$-bundle $\pi:\wt M\lra M$ associated with $E$.
The following assertion is routine.
\begin{proposition}
The space $\wt M$ is a connected manifold with corners equipped with the
natural action of $\cG$ given in any chart $U\times\cG$ on $\wt M$, where
$U$ is a chart on $M$, by the formula $\si(z,g)=(z,g\si^{-1})$,
$\si\in\cG$. The lift $\pi^*E$ is a trivial bundle, $\pi^*E\simeq\wt
M\times\RR^s$, where the trivialization is uniquely determined up to an
automorphism of $\RR^s$. The natural projection $\hat\pi:\pi^*E\lra E$ is
given in coordinates by the formula
\begin{equation*}
    U\times\cG\times\RR^s\ni(z,g,y)\longmapsto (z,gy)\in U\times\RR^s.
\end{equation*}
\end{proposition}

The space $L^2(E)$ (where the measure on $E$ is locally chosen as the
direct product of the measure on $M$ constructed in the preceding
subsection by the standard Lebesgue measure in the fibers) can be
identified with the subspace $L^2_\cG(\pi^*E)\subset L^2(\pi^*E)$ formed by
$\cG$-invariant functions $u(x,y)$, i.e., functions satisfying the
condition
\begin{equation*}
    u(\si x,\si y)=u(x,y), \qquad x\in\wt M,\quad y\in\RR^s,\quad \si\in\cG.
\end{equation*}
\begin{definition}\label{transla}
A bounded operator
$$
 A:L^2(E)\lra L^2(E)
$$
is said to be \textit{translation invariant} if it is the restriction to
$L^2_\cG(\pi^*E)\simeq L^2(E)$ of a bounded $\cG$-invariant operator
\begin{equation}\label{transinv}
 \wt A:L^2(\pi^*E)\lra L^2(\pi^*E)
\end{equation}
that is invariant under translations in
$\RR^s$:
\begin{equation*}
 [\wt Au](x,y+t)=\widetilde{A}[u(x,y+t)],\quad \forall t\in\RR^s.
\end{equation*}
\end{definition}

\begin{proposition}\label{prop-1}
If $A:L^2(E)\lra L^2(E)$ is a translation-invariant operator, then the
corresponding operator~\eqref{transinv} is unique.
\end{proposition}
\begin{proof}
The proof will be given in the appendix.
\end{proof}
\begin{remark}
Since the trivialization of $\pi^*E$ is uniquely determined up to an
automorphism of $\RR^s$ (independent of the point of the base $\wt M$), the
notion of a translation-invariant operator in $\wt E$ is well defined
(automorphisms of $\RR^s$ take translations to translations). It is here
where the requirement of minimality of the structure group is important:
without it, there would be several trivializations of $\pi^*E$ not taken to
each other by a constant automorphism of the fiber, and the notion of
translation-invariant operator would be ambiguous.
\end{remark}

The translation-invariant operator~\eqref{transinv} can be represented in
the form
\begin{equation}\label{sym-tra-inv1}
    \wt A=B\BL(-i\pd{}y\BR),
\end{equation}
where
\begin{equation}\label{sym-tra-inv2}
    B(q):L^2(\wt M)\lra L^2(\wt M),\quad q\in\RR^s,
\end{equation}
is a bounded operator-valued function that is (at least) strongly
measurable in $q$ \cite[Proposition~16]{NaSaSt2}.
\begin{definition}\label{tio}
The function~\eqref{sym-tra-inv2} is called the \textit{symbol} of the
translation-invariant operator $A$ and is denoted by $\si(A)$.
\end{definition}

\subsection{General local operators and localization principle\spacefactor1001}\label{ss22}

\subsubsection*{Local operators with parameters}

Let $X$ be a separable locally compact metric space equipped with a
nonatomic Borel measure $\mu$ such that $\mu(U)>0$ for any nonempty open
set $U\subset X$. We deal with \textit{local operators with a parameter
$q\in\RR^s$} in the $C(X)$-module $H=L^2(X,d\mu)$. They are defined as
operator families $A\in C(\RR^s,\cB H)$ such that for each $\ph\in C_0(X)$
the commutator $[A(q),\ph]$ belongs to the ideal $\cJ=C_0(\RR^s,\cK H)$ of
compact-valued families decaying in norm as $q\to\infty$. Such families $A$
obviously form a $C^*$-subalgebra in $C(\RR^s,\cB H)$, which will be
denoted by $\cA=\cA(\RR^s,\cB H)$.

\subsubsection*{Localization principle}

For $x\in X$, let $\cJ_x\subset\cA$ be the ideal in $\cA$ generated by the
maximal ideal $\cI_x\subset C(X)$ of functions vanishing at $x$, and let
$p_x:\cA\lra\cA_x$ be the natural projection into the \textit{local
algebra} $\cA_x=\cA\slash\cJ_x$.

\begin{theorem}[localization principle; cf.~{\cite[Theorem~3]{NaSaSt2}}]\label{loc-in-cA}
Suppose that the space $X$ is compact. Then $\cJ=\bigcap_{x\in X}\cJ_x$,
and hence an operator $A\in\cA$ is
\begin{enumerate}
    \item Compact with parameter $q$ \rom($A\in\cJ$\rom)
    if and only if all its local representatives
    $p_x(A)\in\cA_x$ are zero.
    \item Fredholm with parameter $q$ \rom(invertible modulo $\cJ$\rom)
    if and only if all its local representatives
    $p_x(A)\in\cA_x$ are invertible.
\end{enumerate}
\end{theorem}
The ideals $\cJ_x$ can be described as follows. For $U\subset X$ and
$A\in\cA$, set
\begin{equation}\label{par-norma}
    \norm{A}_U=\sup_{q\in \RR^s}\norm{A(q)\mid_{H_U}:H_U\lra H},
    \quad\text{where }H_U=\{v\in H\colon \supp v\subset \ov U\}.
\end{equation}
\begin{proposition}[cf.~{\cite[Proposition~4]{NaSaSt2}}]\label{Jx}
The ideal $\cJ_x$ is the set of elements $A\in\cA$ such that
\begin{equation}\label{usl-1}
    \lim_{U\downarrow x}\norm{A}_U=0.
\end{equation}
\rom(Here the limit is taken over the filter of neighborhoods of $x$, i.e.,
over a sequence of open sets $U$ shrinking to $x$.\rom)
\end{proposition}

\begin{remark}
Condition~\eqref{usl-1} is stated in~\cite{NaSaSt2} in the different form
$\lim\norm{A\ph}=0$, where $\abs{\ph}\le1$ and the support of $\ph$ shrinks
to $x$; the two forms are easily seen to be equivalent. The assumption that
$X$ is compact is also easily removed.
\end{remark}

\subsubsection*{Local representatives}

Let us describe the range of the family $\{p_x\}_{x\in X}$ of ``localizing
homomorphisms." Consider a family $\{a_x\in\cA_x\}_{x\in X}$. For each $x$,
we arbitrarily pick up some representative $A_x\in a_x$.
Proposition~\ref{Jx} has an immediate corollary:
\begin{corollary}\label{Jxc}
The family $\{a_x\}$ has the form $a_x=p_x(A)$ for some $A\in\cA$ if and
only if for any $\e>0$ each point $x\in X$ has a neighborhood $U(\e,x)$
such that
\begin{equation}
 \norm{A_x-A}_{U(\e,x)}\le\e.\label{e-usl0}
\end{equation}
\end{corollary}
This is not especially useful, because one has to know $A$ in advance.
Fortunately, one can give a criterion that does not resort to $A$.
\begin{definition}\label{de1}
The family $\{a_x\}$ is said to be \textit{continuous} if for all $\e>0$
and $x\in X$ there exist neighborhoods $U(\e,x)$ such that
\begin{equation}\label{e-usl}
    \norm{A_y-A_{y'}}_{U(\e,y)\cap U(\e,y')}\le\e\quad \text{for any $y,y'\in X$.}
\end{equation}
\end{definition}
One can readily see that the definition of continuity is independent of the
choice of $A_x\in a_x$ (but the neighborhoods $U(\e,x)$ depend on this
choice).
\begin{proposition}[cf.~{\cite[Proposition~7]{NaSaSt2}}]\label{p2}
\rom{(i)} The family $\{a_x\}$ is continuous if \rom(and, in the case of
compact $X$, only if\rom) it has the form $a_x=p_x(A)$ for some $A\in\cA$.

\rom{(ii)} Under the assumptions of \rom{(i)}, if $a_x\in\cB\slash\cJ$ for
all $x\in X$, where $\cB\subset\cA$ is a $C^*$-subalgebra containing $\cJ$,
then $A\in\cB$.
\end{proposition}

\begin{remark*}
For the general localization principle, the topology on the disjoint union
$\bigsqcup_x\cA_x$ in which the families $\{p_x(A)\}_{x\in X}$, $A\in\cA$,
are exactly continuous sections of the projection $\bigsqcup_x\cA_x\lra X$
is described e.g., in~\cite{DaHo1,Vas3}. In our special case, these
sections admit the simpler description given above.
\end{remark*}

\subsubsection*{Infinitesimal operators}
The study of local representatives of an operator $A\in\cA$ is also local
in the following sense. The class $p_{x_0}(A)\in\cA_{x_0}$ remains
unchanged if we multiply $A$ (on the left or on the right) by any cutoff
function $f\in C_0(X)$ such that $f(x_0)=1$. (This can readily be derived
from the fact that if $K\in\cJ$, then $\norm{K}_U\to 0$ as $U\downarrow
x$.) It follows that only what happens in an arbitrarily small neighborhood
of $x_0$ is actually important. Consequently, if $\wt X$ is another metric
space equipped with a measure $\wt\mu$ and a homeomorphism $f\colon
U\lra\wt U$ of some neighborhood $U\subset X$ of $x_0$ onto a neighborhood
$\wt U\subset \wt X$ of the point $\wt x_0=f(x_0)$ is given such that $f$
respects the classes of the measures $\mu$ and $\wt\mu$, then $f^*$ induces
an isomorphism between $\cA_{x_0}$ and $\wt\cA_{\wt x_0}$ and one can speak
of local representatives of $A$ in the algebra $\wt A$ of local operators
with a parameter on $\wt X$.

We systematically use this construction in what follows; the space $\wt X$
will only reflect local properties of $X$ near $x_0$ and is usually
noncompact. Such local representatives, uniquely determined by certain
additional conditions, will also be called \textit{infinitesimal operators}
to emphasize the fact that $X\ne\wt X$.

\begin{example}
If $A$ is a pseudodifferential operator on a smooth manifold $X$, then one
can identify a small neighborhood of $x_0$ with a small neighborhood of
zero in $\wt X=T_{x_0}X$ via the geodesic exponential mapping and take the
operator $\si(A)(x_0,-i\pa/\pa y)$ with constant coefficients on $T_{x_0}X$
for a local representative (infinitesimal operator) of $A$ at $x_0$. (Here
$\si(A)$ is the principal symbol of $A$ and $y\in T_{x_0}X$.) This
infinitesimal operator is uniquely determined by the condition of
invariance with respect to the dilations $y\longmapsto\la y$ in $T_{x_0}X$.
\end{example}

\subsection{Definition and Properties of $\Psi$DO}\label{ss23}

Now we are in a position to define pseudodifferential operators with a
parameter $q\in\RR^s$ on a manifold $M$ with corners. They will be local
operators with a parameter in the sense of Sec.~\ref{ss22} possessing a
number of additional properties.

\subsubsection*{Parameter dependence of $\Psi$DO}

First of all, we impose more restrictive conditions on the dependence of
operators on the parameter $q\in\RR^s$.

We treat $L^2(M)$ as a module over $C(M^\#)$ (by interpreting elements
$u\in L^2(M)$ as functions on $\ciM=M^{\#\circ}$) and consider the algebra
$\cA(\RR^s,\cB L^2(M))$ of operators with parameter $q\in\RR^s$ local with
respect to the action of $C(M^\#)$.
\begin{definition}
The \textit{subalgebra $\cA_{scv}\equiv\cA_{scv}(M,\RR^s)\subset
\cA(\RR^s,\cB L^2(M))$ of functions of slow compact variation} consists of
operator families $A(q)$, $q\in\RR^s$, satisfying the following conditions:
\begin{enumerate}
    \item The function $A(q)$ is of \textit{compact variation}; that is,
    $A(q)-A(q')\in\cK H$ for any $q,q'\in\RR^s$.
    \item The function $A(q)$ \textit{slowly varies at infinity}
    in the sense that for any $d>0$ and $\e>0$ there exists an $R>0$ such that
\begin{equation*}
    \norm{A(y)-A(y')}\le\e \quad\text{whenever $\abs{y-y'}<d$ and $\abs{y}>R$}.
\end{equation*}
\end{enumerate}
\end{definition}
\begin{proposition}\label{cAsv}
The set $\cA_{scv}$ is a $C^*$-algebra, and every element $A(y)\in\cB$ can
be approximated by $C^\infty$ functions of compact variation all of whose
derivatives decay at infinity.
\end{proposition}
\begin{proof}
The proof will be given in the appendix.
\end{proof}

\subsubsection*{Interior symbol}

Let $A\in\cA_{scv}(M,\RR^s)$.
\begin{definition}\label{AV}
Let $x\in \ciM$ be an interior point of $M$. We say that $A$ is
\textit{Agranovich--Vishik at $x$} if, under the identification of a
neighborhood of $x$ in $\ciM$ with a neighborhood of the origin in $T_{x}M$
via a coordinate system near $x_0$, $A$ has a local representative of the
form
\begin{equation*}
    A_{x_0}=B\BL(q,-i\pd{}y\BR),\quad y\in T_{x_0}M,
\end{equation*}
where $B(q,\xi)$ is a function continuous for $\abs{q}+\abs{\xi}\ne0$ and
zero-order homogeneous:
\begin{equation*}
    B(\la q,\la \xi)=B(q,\xi),\quad \la\in\RR_+.
\end{equation*}
The function $B(q,\xi)$ is called the \textit{interior symbol} of $A$ and
is denoted by
\begin{equation*}
    \si_0(A)(x,\xi,q):=B(q,\xi).
\end{equation*}
\end{definition}
Essentially, the definition says that at the point $x$ the operator $A$ is
a parameter-dependent pseudodifferential operator with continuous symbol.
\begin{proposition}
If $A$ is Agranovich--Vishik at $x$, then $\si_0(A)$ is a well-defined
function on $T^*_xM\times\RR^s$ outside zero (i.e., its existence and form
is independent of the choice of the coordinate system).
\end{proposition}
\begin{proof}[Sketch of proof]
The operator
\begin{equation*}
 \wh B=B\BL(q,-i\pd{}y\BR)
\end{equation*}
behaves as desired under linear changes of the variable $y$. Thus,
essentially, one should prove that if $f:\RR^n\to\RR^n$ is a diffeomorphism
with identity differential at the origin, then $\wh B$ and $(f^*)^{-1}\wh
Bf^*$ define the same element in the local algebra $\cA_0$. To prove this,
we approximate $B$ by smooth classical symbols and use the theorem on the
change of variables in a classical pseudodifferential operator.
\end{proof}

\subsubsection*{Face symbols}

From now on, we choose and fix a compatible system of exponential maps from
the normal bundles of the faces into their neighborhoods in $M$. Our
definition of face symbols and of $\psi$DO tacitly depends on the choice of
this system. The questions concerning the invariance of the definition will
be discussed elsewhere.

Let again $A\in\cA_{scv}(M,\RR^s)$, and let $z\in F^\#$ be a point of open
face $F^\#$ dual to a face $F$ of positive codimension $d\ge 1$ in $M$.
Some neighborhood $U$ of $F$ can be identified via the exponential map with
a neighborhood of the zero section in $N_+F$ or (additionally applying the
logarithmic map) with a neighborhood of the point at infinity in the
positive quadrant in$NF$. Hence we have the embedding
$$
L^2(M)|_U\subset L^2(N_+F)\simeq L^2(NF),
$$
which implies that it suffices to localize the operator A at the point
$z\in F^\#$ in the space $L^2(NF)$.
\begin{definition}\label{face-sym}
We say that the operator $A(q)$ has a \textit{translation-invariant
infinitesimal operator} at the point $z\in F^\#$ if in $L^2(NF)$ there
exists a translation-invariant operator $A_\infty(q)$ (see
Definition~\ref{transla}) belonging to the same coset in the local algebra
$\cA_z$. The symbol (Definition~\ref{tio}) of $A_\infty(q)$ will be called
the \textit{symbol of $A(q)$ at} $z$ and will be denoted by $\si_z(A)$.
\end{definition}

\begin{theorem}\label{grasy}
If $A(q)$ has a translation-invariant infinitesimal operator, then it is
unique. Thus the symbol $\si_z(A)$ is well defined. It is a
$\gS_F$-invariant operator-valued function on $\RR^d\times\RR^s$ with
values in $\cB L^2(\wt F)$, where $\wt F$ is the principal $\gS_F$-covering
over $F$ trivializing $NF$.
\end{theorem}

\begin{proof}
[Proof \/\rm will be given in the Appendix]
\end{proof}

\subsubsection*{Pseudodifferential operators}

Let $M$ be a manifold with corners.

\begin{definition}\label{P-D-O}
The space $\Psi(M)\equiv\Psi(M,\RR^s)$ of pseudodifferential operators
consists of operator families $A(q)$ satisfying the following
conditions:
\begin{enumerate}
    \item $A(q)\in\cA_{scv}(M,\RR^s)$.
    \item For each interior point $x\in M^\circ$, the family $A(q)$ is
    Agranovich--Vishik at $x$.
    \item For each face $F$ of codimension $d=d(F)>0$ in $M$, the family
    $A(q)$ has a $\gS_F$-invariant symbol $\si_z(A)$ in the sense of
    Definition~\ref{face-sym} at each point $z\in F^\#$, and $\sigma_z(A)$ is
    independent of $z$. Moreover,
    $\si_z(A)\in\Psi(\wt F,\RR^{d+s}$); i.e., the symbol $\si_z(A)$
    is a $\gS_F$-invariant $\psi$DO with parameters $(q,p)\in\RR^s\times\RR^d$
    on the manifold
    $\wt
    {F}$ with corners, the covering of $\ov{F}$ trivializing the bundle $N\ov{F}$.
\end{enumerate}
Since the symbol $\si_z(A)$ is independent of $z\in F^\#$, it will be
denoted by $\si_F(A)$ in what follows. The interior symbol will be denoted
by $\si_0(A)$; it is defined on the interior of $T^*M\times\RR^s$ minus the
zero section.
\end{definition}

\subsubsection*{Main theorem of the calculus}

The localization principle (Theorem~\ref{loc-in-cA}) readily implies the
following assertion.
\begin{theorem}[main theorem of the calculus]
A pseudodifferential operator $A$ on a compact manifold $M$ with corners is
uniquely determined modulo the ideal $\cJ$ of compact operators with
parameters by the symbol tuple $(\si_0(A),\{\si_F(A)\})$, where $F$ runs
over all faces of positive codimension, modulo compact operators. The map
\begin{equation*}
 \si\colon A\longmapsto(\si_0(A),\{\si_F(A)\})
\end{equation*}
that takes each $\psi$DO $A\in\Psi(M,\RR^s)$ to it symbol tuple is a
$C^*$-algebra homomorphism.
\end{theorem}

\subsubsection*{The symbol algebra}

Now let us describe the symbol algebra, i.e., the range of the symbol map
$\si$. In other words, we should indicate conditions on the interior symbol
and the face symbols on faces of positive codimension necessary and
sufficient for the existence of a $\psi$DO with
these symbols. To avoid awkward formulas, we first do so for the case in
which the normal bundles of all faces are trivial and then indicate the
modifications needed in the general case.

Thus let $M$ be a manifold with corners such that the normal bundle $NF$ is
trivial for all faces $F$ of $M$.

Let the following data be given:
\begin{itemize}
    \item
For each interior point $x\in\ciM$, a continuous zero-order homogeneous
function $\si_x$ on $(T^*_xM\times\RR^s)\setminus0$.
    \item
For each face $F$ of codimension $d>0$, a pseudodifferential operator
$\si_F\in\Psi(F,\RR^{d+s})$.
\end{itemize}

\begin{theorem}[description of the symbol algebra]
For the existence of a $\psi$DO $A\in\Psi(M,\RR^s)$ such that
\begin{align}
    \si_0(A)&=\si_x \quad\text{on $(T^*_xM\times\RR^s)\setminus0$ for each $x\in\ciM$,}
    \label{lab1}\\
    \si_F(A)&=\si_F \quad\text{for each face $F$ of positive codimension,}
    \label{lab2}
\end{align}
the following conditions are necessary and sufficient\rom:
\begin{enumerate}
    \item
The functions $\si_x$ form a continuous zero-order homogeneous function on
the interior of $(T^*M\times\RR^s)\setminus0$ and extend by continuity to a
continuous function \rom(which we denote by $\si_0$\rom) on the whole space
$(T^*M\times\RR^s)\setminus0$.
    \item The restriction of $\si_0$ to the boundary satisfies the
    compatibility conditions
\begin{equation}\label{comp1}
    \si_0\bigm|_{F}=\si_0(\si_F)  \quad\text{for each face $F$ of positive codimension,}
\end{equation}
where the left-hand side is the restriction of $\si_0$ to
$T^*M|_F\oplus\RR^s$, naturally identified with $T^*F\oplus
N^*F\oplus\RR^s=T^*F\oplus\RR^{d+s}$.
    \item
If $F_1\succ F_2$ are two adjacent faces of $M$ and $\Ga$ is a face of
$F_1$ mapped into $F_2$ under the immersion of $F_1$ in $M$, then
\begin{equation}\label{comp2}
    \si_\Ga(\si_{F_1})=\si_{F_2}.
\end{equation}
\end{enumerate}
\end{theorem}
\begin{proof}
First, note that routine computations based on composition formulas for
pseudodifferential operators and standard norm estimates show that, being
quantized, the symbols $\si_{x_0}(q,p)$ and $\si_F(q,\xi)$ give rise to the
local representatives $\wh\si_x=\si_{x_0}(q,-i\pa/\pa x)$ and
$\wh\si_F=\si_F(q,-i\pa/\pa t)$ that  belong to $\cA_{csv}$.

By Proposition~\ref{p2}, to prove the theorem it remains to establish that
conditions (1)--(3) are exactly equivalent to the continuity of this family
of local representatives in the sense of Definition~\ref{de1}.

(a) Let us show that the function $\si_x$ continuously depends on $x$ in
the interior of $M$. Localizing our considerations, we can assume that
$M=\RR^n$. The family $\si_x$ is continuous if and only if for each $\e>0$
each point has a neighborhood $U(\e,x)$ such that
$\norm{\wh\si_x-\wh\si_y}_{U(\e,x)\cap U(\e,y)}\le\e$ for any $x$ and $y$.
The intersection $U=U(\e,x)\cap U(\e,y)$ is necessarily nonempty if $y\in
U(\e,x)$. Since the operator $\wh\si_x-\wh\si_y$ is dilatation invariant,
it follows that
\begin{equation*}
 \norm{\wh\si_x-\wh\si_y}_U=\norm{\wh\si_x-\wh\si_y}=\max_p\abs{\si_x-\si_y}
\end{equation*}
(provided that $U$ is nonempty). Combining this with the homogeneity of
$\si_x$ in $(p,q)$, we see that the continuity of the family of local
representatives in the sense of Definition~\ref{de1} is equivalent to the
continuity of the interior symbol. This is of course well known from the
theory of $\psi$DO on smooth manifolds.

(b) Let us show that the interior symbol is continuous up to the boundary
and satisfies the compatibility conditions~\eqref{comp1} there. Fix a point
$z_0\in F$. Multiplying by a cutoff function $f\in C(M)$, we can study the
problem assuming that $F=\RR^{n-d}$ and $\ciM=\RR^{n-d}\times\RR^d$ (here
we use the logarithmic coordinates $y\in\mathbb{R}^d$, see \eqref{loga1} on
the fibers of the normal bundle of $F$). Let $x_0\in F$. Applying
Corollary~\ref{Jxc} and using the Fourier transform with respect to
$\RR^d$, we see that for each $\e>0$ there is a neighborhood $U_\e$ of
$x_0$ in $F$ such that
\begin{equation}\label{a1}
    \norm{\wh\si_F-\wh\si_{x_0}(\si_F)}_{U_\e\times\RR^d}<\e.
\end{equation}
(Here $\si_{x_0}(\si_F)$ is the symbol $\si_0(\si_F)$ restricted to the
fiber over $x_0$.)

On the other hand, the continuity of the family of local representatives on
$M$ near $F^\#$ is equivalent to the existence of a neighborhood $W_\e$ of
the point at infinity on the diagonal of the positive quadrant in $\RR^d$
such that
\begin{equation}\label{a2}
    \norm{\wh\si_F-\wh\si_{x}}_{(F\times W_\e)\cap U(\e,x)}<\e.
\end{equation}
For $x\in U_\e\times W_e$, using the triangle inequality, from~\eqref{a1}
and \eqref{a2} we conclude that
\begin{equation}\label{a4}
    \norm{\wh\si_x-\wh\si_{x_0}(\si_F)}_{U}<2\e
\end{equation}
on the nonempty set $U=(U_\e\times W_e)\cap U(\e,x)$. Arguing as above, we
see that $\abs{\si_x-\si_{x_0}(\si_F)}\le2\e$ for these $x$.

(c) In a similar way, one shows that condition~(3) also follows from the
continuity of local representatives and finally concludes that
conditions~(1)--(3) together are equivalent to the continuity. We leave the
details to the reader.
\end{proof}

\begin{remark}
In particular, it follows from the compatibility condition that the symbol
on a face of positive codimension determines the symbols on all adjacent
faces of larger codimension.
\end{remark}

Let us now discuss how the compatibility conditions should be modified if
the normal bundles of the faces are not trivial.

Let again $\ov{F}_1\succ \ov{F}_2$ be two adjacent faces of $M$, and let
$\Ga$ be a face of $\ov{F}_1$ covering $\ov{F}_2$ (there can be several
such faces). The symbols $\si_{F_1}(A)$ and $\si_{F_2}(A)$ of $A$ are
operators with parameters on the minimum coverings $\wt F_1$ and $\wt F_2$
trivializing the bundles $NF_1$ and $NF_2$, respectively. Let $\wt\Ga$ be
the lift of $\Ga$ to $\wt{F}_1$. The symbol $\si_{\wt\Ga}(\si_{F_1}(A))$ is
defined on the covering $\wt{\wt\Ga}$ trivializing the bundle $N\wt\Ga$.
The composite covering $\wt{\wt\Ga}\lra F_2$ trivializes $NF_2$ (since it
trivializes both direct summands $N\wt\Ga$ and $NF_1|_\Ga$). Since the
trivializing covering $\wt F_2\lra F_2$ is minimal and hence universal,
there exists a unique (modulo permutation of the sheets) covering
$\wt{\wt\Ga}-\to\wt F_2$ making the triangle
\begin{equation}\label{treug}
     \xymatrix{
  \wt{\wt\Ga} \ar[dr] \ar@{-->}[r]^\pi
                & \wt F_2 \ar[d]   \\
                & F_2             }
\end{equation}
commute.

Let $L^2_{inv}(\pi)$ be the subspace of $L^2(\wt\Ga)$ consisting of
functions invariant with respect to permutations of the sheets of $\pi$.
The compatibility condition~\eqref{comp2} in this situation is generalized
to
\begin{equation}\label{soglas1}
    \si_{\wt\Ga}(\si_{F_1}(A))|_{L^2_{inv}(\pi)}=\si_{F_2}(A).
\end{equation}

The counterpart of the compatibility condition~\eqref{comp1} reads
\begin{equation}\label{soglas2}
    \si_0(\si_F(A))=\pi_F^*\left[\si_0(A)|_{T^*F}\right],
\end{equation}
where $\pi_F:T^*\wt F\lra T^*F$ is the covering associated with the
covering $\wt F\lra F$.

\appendix

\section{Proofs of Some Assertions}

\subsubsection*{Proof of Proposition~\rom{\ref{p-1}}}

For brevity, we write
\begin{equation*}
  F=\cGa_j(M)\quad
d=d_j.
\end{equation*}

Let $U\simeq \ov\RR{}_+^s\times\RR^{n-s}$ be a coordinate neighborhood on
$M$. If the intersection $U\cap F$ is nonempty (which can happen only for
$s\ge d$), then it consists of finitely many ($\le C_s^d$) connected
components of the form $V\simeq \RR_+^{s-d}\times\RR^{n-s}$, where the open
coordinate quadrant $\RR_+^{s-d}$ of dimension $s-d$ is singled out in
$\ov\RR{}_+^s$ by the relations
\begin{equation*}
 x_{j_{1}}=\dotsm=x_{j_d}=0,\quad x_{j_{d+1}},\dotsc,x_{j_{s}}>0
\end{equation*}
for some (depending on $V$) permutation $j_1,\dotsc,j_s$ of the indices
$1,\dotsc,s$. If we accordingly permute the standard coordinates
$x_1,\dotsc x_n$ in $U$, setting
\begin{multline*}
 \rho_1=x_{j_{1}},\dotsc,\rho_d=x_{j_{d}},\\
y_{d+1}=x_{j_{d+1}},\dotsc,y_s=x_{j_{s}}, y_{s+1}=x_{s+1}, \dotsc, y_n=x_n,
\end{multline*}
then the variables $y=(y_{d+1},\dotsc,y_n)$ are coordinates in $V$ and the
variables $\rho=(\rho_1,\dotsc,\rho_d)$ are defining functions of $V$
for the embedding $V\subset U$ (and local defining functions of
$\ov{F}$); i.e., locally the face is given by the conditions $\rho=0.$

We take a finite cover of $M$ by coordinate neighborhoods $U$ and various
connected components $V\subset U\cap F$ and obtain a finite atlas
\begin{equation*}
 \bl\{\bl(V,\,y:V\lra\RR{}_+^{s-d}\times\RR^{n-s}\br)\br\}
\end{equation*}
on $\ov{F}$ such that associated with each coordinate neighborhood $V$ of
this atlas is a pair $(U,V)$ and coordinates $(\rho,y)$ in $U$. Let $\wt V$
be another coordinate neighborhood $(\wt U,\wt V)$ with coordinates $(\wt
\rho,\wt y)$ in $\wt U$, and suppose that the intersection $V\cap\wt V$ is
nonempty. The change of variables
\begin{equation*}
    \wt y\circ y^{-1}:y(V\cap\wt V)\lra\wt y(V\cap\wt V)
\end{equation*}
is obtained by restriction to $y(V\cap\wt V)$ from the change of
coordinates $(\rho,y)\longmapsto(\wt\rho,\wt y)$ on the intersection of the
coordinate neighborhoods $U$ and $\wt U$ on $M$ and hence has a smooth
continuation to the closure of the set $y(V\cap\wt V)$ in
$\ov\RR{}_+^{s-d}\times\RR^{n-s}$. (The continuation is obtained by
restriction of the same change of coordinates to the closure.) These
continuations determine the transition functions of some compact manifold
${\ov{F}}$ with corners whose local models are
$\ov\RR{}_+^{s-d}\times\RR^{n-s}$ and into which $ F$ is naturally embedded
as a dense open submanifold. This manifold $\ov{F}=:\Ga_j(M)$ is the
\textit{closed face} of $M$ corresponding to the open face $\cGa_j(M)$. The
embedding $\cGa_j(M)\subset M$ extends by continuity to $\Ga_j(M)$; the
resulting mapping is in general an immersion with self-intersections.

\subsubsection*{Proof of Proposition~\rom{\ref{hahaha}}} It suffices to write out
a natural invariant pairing; this can be done in the coordinates
$(\rho,y)$: for a form
$$
\om=\sum a_j\rho_j^{-1}d\rho_j\in N^*\ov{F}
$$
and a vector
$$
\xi=(b_1,\dotsc,b_d)\in N\ov{F},
$$
we set
$$
\langle\om,\xi\rangle=\sum a_jb_j.
$$
Under changes of coordinates, the components of $\xi$ and $\om$ are subjected
to the same permutation, and the defining functions $\rho_j$ are multiplied
by nonzero numbers (the diagonal entries of the matrix $\La(y)$), which does
not affect the logarithmic derivatives, so that the numbers $a_j$ remain the
same. Thus the pairing is independent of the choice of coordinates.

\subsubsection*{Proof of Theorem~\rom{\ref{papkahoroshaja}}}

We need the following simple lemma.
\begin{lemma}
If smooth mappings
$$
g_j:\RR_+^k\lra\RR_+^k,\quad g_j(0)=0,\, j=1,\dotsc,l,
$$
are diffeomorphisms in a neighborhood of the origin, all matrices
$g_j'(0)\bl(g_i'(0)\br)^{-1}$ are diagonal, and $\la_1,\dotsc,\la_l$ are
nonnegative numbers at least one of which is nonzero, then the mapping
\begin{equation*}
    g\equiv\sum_{j=1}^l\la_jg_j:\RR_+^k\lra\RR_+^k
\end{equation*}
is also a diffeomorphism in a neighborhood of the origin.
\end{lemma}
Indeed, the only nontrivial assertion is that $g$ is epimorphic, but this
can be verified as follows. Since the matrices
$g_j'(0)\bl(g_i'(0)\br)^{-1}$ are diagonal, it follows that all $g_j$ take
any given coordinate quadrant of arbitrary dimension to one and the same
coordinate quadrant.

\medskip

The lemma suggests that one can construct the desired mapping $f=f_j$
specifying it locally by the formula
\begin{equation}\label{bezim}
\rho=r,
\end{equation}
where $\rho$ is a local tuple of defining functions of the face
$\Gamma_j(M)$ and $r$ are the corresponding coordinates in the fiber of
$N_+\Ga_j(M)$ and then gluing the local mappings with the use of a partition
of unity.

We implement this idea and construct the mapping $f$, successively
extending the set on which it is defined. Suppose that $f$ has already been
defined over some set $O\subset\Ga_j(M)$, and let $V$ be a local chart on
$\Ga_j(M)$ with the corresponding pair
$(U=V\times\overline{\mathbb{R}}^{d_j}_+,V)$, so that over $V$ the mapping
can be given by formula \eqref{bezim}. Let $(\ph_O,\ph_V)$ be a nonnegative
partition of unity on $O\cup V$ subject to the cover by $O$ and $V$. we
construct the map over $O\cup V$ by setting
\begin{equation*}
    f_{O\cup V}=\begin{cases}
    f_O &\text{over $O\setminus\supp\ph_V$,}\\
    \ph_Of_O +\ph_Vf_V &\text{over $V$,}
    \end{cases}
\end{equation*}
where the addition in the second line is carried out in the fibers of $U\lra V$
(and is well defined in a sufficiently small neighborhood of $V$). Since
$\ph_O=1$ on $V\setminus\supp\ph_V$, it follows that both definitions are
compatible on the set where they apply simultaneously, and the lemma now
implies that we have defined a mapping with the desired properties over $O\cup
V$.

To complete the proof of Theorem~\ref{papkahoroshaja}, it suffices to start
from an empty set $O$ and successively add to it all charts from a finite
atlas on $\Ga_j(M)$. To obtain compatible (in the sense that diagram
\eqref{dia1} commutes) exponential mappings for all faces, one should start
from faces of maximal codimension.

\subsubsection*{Proof of Proposition~\rom{\ref{prop-1}}} Let us proceed from the
operator $\wt A$ to the symbol $B(p)=\si(A)$. We have to prove its
uniqueness; it is a consequence of the following assertion, which we state
in general form.

\begin{lemma}\label{eq2}
Let a finite group $G$ act on the space $L^2(\RR^s;H)$, where $H$ is a
Hilbert space, by the formula
\begin{equation*}
    [T_gf](p)=S_gf(\si_g^{-1}p),\quad p\in\RR^s,
\end{equation*}
where $S$ is a representation of $G$ on $H$ and $\si$ is a
\textbf{faithful} representation of the same group on $\RR^s_p$. Let
\begin{equation*}
    A(p)\colon H\lra H,\quad p\in\RR^s,
\end{equation*}
be a continuous operator-valued function. Then the following assertions
hold.
\begin{enumerate}
 \item[\rom{(ii)}]  If $A(p)f(p)=0$ for almost all $p$ for any element $f\in
L^2(\RR^s;H)$ such that
\begin{equation}\label{muka1}
    T_gf=f \quad\forall g\in G,
\end{equation}
then $A(p)=0$ for all $p$.
  \item[\rom{(ii)}] If the operator
$A:L^2(\RR^s;H)\longrightarrow L^2(\RR^s;H)$ induced by the multiplication
by $A(p)$ preserves the subspace of invariant functions \eqref{muka1}, then
it is $G$-invariant, i.e., satisfies
$$
A(p)=S_g^{-1}A(\sigma_g(p))S_g,\qquad \forall g\in G.
$$
\end{enumerate}
\end{lemma}
\begin{proof}
1. Since $A(p)$ is continuous, it suffices to prove the desired relation on
a dense set of values of $p$. For this set we take $\Om=
\RR^s\setminus\bl(\bigcup_j\op{fix}\si_g\br)$, where $\op{fix}\si_g$ is the
set of fixed points of $\si_g$ (which is at most a hyperplane, since $\si$
is faithful). Let $p_0\in\Om$ and $v\in H$. We claim that $A(p_0)v=0$.
Indeed, let $U$ be a sufficiently small neighborhood of $p_0$ such that
$U\subset\Om$ and hence
\begin{equation*}
    \si_g(U)\cap\si_h(U)=\varnothing\quad\text{for}\quad g\ne h,\quad g,h\in
    G.
\end{equation*}
Then the function
\begin{equation*}
    f(p)=\begin{cases}
    S_gv & \text{if } p\in\si_g(U)\quad\text{for some $g\in G$},\\
    0&\text{otherwise}
    \end{cases}
\end{equation*}
is well defined. It is $G$-invariant, so that $A(p)f(p)=0$ for almost all
$p$ and hence $A(p_0)v=0$ (since $A(p)f(p)\equiv A(p)v$ is continuous in
$U$).

2. The second assertion of the lemma is proved by the same method. Namely,
it follows from the $G$-invariance of $f(p)$ and the assumptions of the
lemma that $T_h(A(p)f(p))=A(p)f(p)$ for each $h\in G$. In turn, this
implies $T_hA(p)T_h^{-1}f(p)=A(p)f(p)$, whence for $p=p_0$ we obtain the
desired relation
$$
S_hA(\sigma_{h^{-1}}(p_0))S_h^{-1}v=A(p_0)v.
$$
\end{proof}
\begin{remark}
It is important that the representation $\si$ is faithful. Without this
assumption, the lemma fails. (One can only prove that $A(p)v=0$ for
elements $v$ invariant under $S_g$ for $g\in\ker\si$). In terms of our
problem, this means that one should always reduce the structure group of th
bundle $N_+\ov{F}$ to the \textit{minimum possible} subgroup.
\end{remark}

\subsubsection*{Proof of Proposition~\rom{\ref{cAsv}}} The first assertion is
trivial. The proof of the second assertion goes by the following scheme:
\begin{enumerate}
    \item The approximation is defined in the standard way as the convolution with a
smooth function with unit integral and small support.
    \item All derivatives of an
approximating function are bounded.
    \item Since the original family slowly varies
at infinity, it follows that the first derivative of the approximating
family decays at infinity.
    \item In turn, (2) and (3) imply that all derivatives decay
at infinity.\qed
\end{enumerate}

\subsubsection*{Proof of Theorem~\rom{\ref{grasy}}}
It suffices to prove uniqueness for the case in which $A(q)\equiv0$.
Consider a sequence $\ph_n\in C_0(N\ov{F})$ strongly convergent to the
identity operator. Then $\ph_n A_\infty\ph_n$ strongly converges to
$A_\infty$. Passing to the cover $N\widetilde{F}$, we see that the product
$\wt\ph_n \wt A_\infty\wt \ph_n$ of the corresponding lifted operators
strongly converges to $\wt A_\infty$. On the other hand, for a given $n$,
let $a_j\in\RR^d$ be a sequence of vectors such that the supports of the
functions $t_{a_j}^*\wt\ph_n$, where $t_{a_j}$ is the shift by the vector
$a_j\in\mathbb{R}^d$, lie as $j\to\infty$ in an arbitrarily small
neighborhood of some of the preimages $z_*$ of the point $z$ (i.e., go to
infinity in the positive quadrant along the corresponding ray). These
functions are no longer $\gS_F$-invariant. However, one can show that there
exist functions $\psi_j$ bounded by $1$ with supports shrinking to $z$ such
that their invariant lifts satisfy the condition
$$
 \wt\psi_jt_{a_j}^*\wt\ph_n=t_{a_j}^*\wt\ph_n.
$$
Then, according to the properties of the local algebra $\mathcal{A}_z$, we
have
\begin{multline*}
(t_{a_j}^*\wt\ph_n)\wt
A_\infty(t_{a_j}^*\wt\ph_n)=(t_{a_j}^*\wt\ph_n)\wt\psi_j\wt
A_\infty\wt\psi_j(t_{a_j}^*\wt\ph_n)
 \\
=(t_{a_j}^*\wt\ph_n)\wt{\psi_j A_\infty \psi_j}(t_{a_j}^*\wt\ph_n)\to0
\end{multline*}
as $j\to\infty$ (convergence in norm). Indeed, the extreme factors are
uniformly bounded, and the middle factor converges to zero, since
$A_\infty$ represents the zero class. Thus
\begin{equation*}
    (t_{a_j}^*\wt\ph_n)\wt
A_\infty(t_{a_j}^*\wt\ph_n)= t_{a_j}^*\circ\wt\ph_n \wt A_\infty
\wt\ph_n\circ t_{a_j}^{*-1}\to 0.
\end{equation*}
(We have used the translation invariance of $A$.) We see that $\wt\ph_n \wt
A_\infty \wt\ph_n=0$ and, passing to the limit as $n\to\infty$, find that
$\wt A_\infty$ and hence $A_\infty$ are zero.

\providecommand{\bysame}{\leavevmode\hbox to3em{\hrulefill}\thinspace}

\end{document}